\documentclass[12pt, reqno, twoside]{amsart}
\usepackage{amsmath, amsthm, amsopn, amssymb, microtype, enumerate}
\usepackage{mathrsfs} 
\usepackage{mathscinet}
\usepackage{bbm}
\usepackage{enumerate, tikz, etoolbox, intcalc, geometry, caption, subcaption, afterpage}
\usepackage[english]{babel}
\usepackage{float}
\usepackage[all]{xy}
\usepackage[pagewise]{lineno}

\usepackage[colorlinks=true,urlcolor=blue,pdfborder={0 0 0}]{hyperref}
\hypersetup{linkcolor=[rgb]{0,0,0.6}}
\hypersetup{citecolor=[rgb]{0,0.6,0}}
\usepackage[
    style=alphabetic,
    giveninits=true,
    sorting=nyt
]{biblatex}

\DeclareFieldFormat[article]{title}{#1}

\renewbibmacro*{in:}{%
  \ifentrytype{article}{}% do nothing
  {\printtext{\bibstring{in}\intitlepunct}}% normal behavior for other types
}

\calclayout

\xyoption{dvips}
\xyoption{color}

\usepackage{color}

\title{An Exceptional Set of Uniformly Spread Kakutani Tilings of the Line}

\author{Yotam Smilansky}
\address{Department of Mathematics,
University of Manchester, United Kingdom.}
\email{yotam.smilansky@manchester.ac.uk}

\newcommand{\PP}{\mathcal{P}}

\newcommand{\TT}{\mathcal{T}}

\newcommand{\N}{{\mathbb{N}}}
\newcommand{\Z}{{\mathbb{Z}}}
\newcommand{\Q}{{\mathbb{Q}}}
\newcommand{\R}{{\mathbb{R}}}
\newcommand{\C}{{\mathbb{C}}}

\newcommand{\absolute}[1] {\left|{#1}\right|}

\newcommand{\disc}[3]{\mathbf{disc}_{#1}\left({#2},{#3}\right)}

\newcommand {\ignore}[1]  {}

\theoremstyle{plain}
\newtheorem{thm}{Theorem}[section]
\newtheorem*{thm*}{Theorem}

\newtheorem{prop}[thm]{Proposition}
\newtheorem*{prop*}{Proposition}

\theoremstyle{definition}
\newtheorem{definition}[thm]{Definition}

\newtheorem{remark}[thm]{Remark}

\newtheorem{example}[thm]{Example}

\newtheorem*{remark*}{Remark}

\numberwithin{equation}{section}
\swapnumbers

\newif\ifdraft\drafttrue
\begin{document}
\maketitle

%\renewcommand{\thefootnote}{}
%\footnotetext{University of Manchester, Department of Mathematics. \texttt{yotam.smilansky@manchester.ac.uk}.}

\begin{abstract}
   The $\alpha$-Kakutani substitution rule splits the unit interval into two subintervals of lengths $\alpha$ and $1-\alpha$, for a fixed $\alpha\in(0,1)$. A simple inflation-substitution procedure produces tilings of the real line and their associated Delone sets. We show that there are precisely five distinct values of $\min\{\alpha,1-\alpha\}$ for which these sets are uniformly spread, meaning that they are a bounded displacement of a lattice. The proof of this surprising fact combines the construction and analysis of a related family of primitive substitution tilings, Solomon's criterion for uniform spreadness, and a classification of Pisot–Vijayaraghavan polynomials. 
\end{abstract}

\section{Introduction}

A set $\Lambda\subset\R^d$ is \emph{Delone} if it is \emph{uniformly discrete} and \emph{relatively dense}, that is, if there exist constants $0<r,R<\infty$ so that $B(x,r)\cap\Lambda=\{x\}$ and $B(y,R)\cap\Lambda\neq\emptyset$ for all $x\in\Lambda$ and $y\in \R^d$, respectively. Clearly, lattices and their translations are Delone sets in $\R^d$, and there are numerous non-lattice examples that are not periodic in any way. Given such a non-lattice Delone set, it is natural to ask how ``close'' it is to a lattice in $\R^d$. One natural interpretation of such a question involves the concept of \emph{bounded displacement (BD) equivalence}, where two Delone sets $\Lambda,\Gamma\subset\R^d$ are BD equivalent if there exists a bijection $\varphi:\Lambda\rightarrow \Gamma$ satisfying $\sup_{x\in\Lambda}\|x-\varphi(x)\|<\infty$. 
A Delone set $\Lambda\subset\R^d$ is \emph{uniformly spread} if it is BD equivalent to some lattice in $\R^d$, or equivalently, by a simple application of Hall's marriage theorem, if there is some $c>0$ for which $\Lambda$ and $c\cdot\Z^d$ are BD equivalent. More generally, two lattices are BD equivalent if and only if they have the same covolume, see \cite{SmilanskySolomon-2025} for a comprehensive discussion.

The Delone sets considered in this manuscript are those associated with $\alpha$-Kakutani tilings for $\alpha\in(0,1)$, a simple class of multiscale substitution tilings \cite{SmilanskySolomon-2021} of $\R$. These tilings are defined by a substitution rule on a closed interval $I\subset\R$ of unit length that splits $I$ into two subintervals, a left one of length $\alpha$ and a right one of length $1-\alpha$. 
 Starting with a copy of $I$, which we think of as a \emph{tile}, we apply a uniform inflation and substitute any tile of length greater than $1$. This results in a growing family of patches that exhausts the real line, and  limits taken in an appropriate way define a space of tilings of $\R$, all made of interval-tiles of lengths between $\min\{\alpha,1-\alpha\}$ and $1$. Given such a tiling, an \emph{associated} Delone set is obtained by selecting one point from each tile, with the requirements that distinct tiles yield distinct points and that the resulting set is uniformly discrete. A natural choice is to take the left endpoint of each tile. However, since all Delone sets associated with a fixed tiling are BD equivalent, the particular choice of points is immaterial. 

Our main result is the following classification.

\begin{thm}\label{thm:main}
    Let $\Lambda_\alpha$ be a Delone set associated with an $\alpha$-Kakutani tiling of $\R$. Then $\Lambda_\alpha$ is uniformly spread if and only if 
    \begin{equation}\label{eq:special values}
            r_\alpha:=\frac{\log(\min\{\alpha,1-\alpha\})}{\log(1-\min\{\alpha,1-\alpha\})}\in\left\{1,\frac{3}{2},2,3,4  \right\}.
    \end{equation}
\end{thm}

For simplicity of presentation, we will assume throughout that $\alpha\in(0,1/2]$, and so $\alpha=\min\{\alpha,1-\alpha\}$. We note that if $r_\alpha=1$ then $\alpha=1/2$, and the left endpoints of tiles in any corresponding tiling are simply translated lattices in $\R$. When $r_\alpha=2$ then $\alpha=1/\phi^2$, where $\phi$ is the \emph{golden ratio}. The other values of $\alpha$ are harder to present explicitly, though the cases $r_\alpha=3$ and $r_\alpha=3/2$ are related to the \emph{supergolden ratio} and the \emph{plastic ratio}, respectively. The latter is the smallest \emph{Pisot–Vijayaraghavan (PV) number}, and the case $r_\alpha=4$ is related to the second smallest PV number, see \S \ref{sec:remarks}. For all other values of $\alpha$ the corresponding Delone sets are never uniformly spread. 

\subsection{Key steps of the proof} 

A complete proof requires some further preparation, but  we can already outline the main steps of the argument. First, in the \emph{incommensurable} case $r_\alpha\notin\Q$, the discrepancy in $\alpha$-Kakutani tilings, that is, the fluctuation of the number of tiles that appear in large patches relative to their expected number (namely, the density), can be shown to be large. This, when combined with Laczkovich's criterion for uniform spreadness \cite{Laczkovich-1992}, is enough to establish that any associated Delone set is not BD equivalent to a lattice. This is a special case of \cite[Theorem 8.2]{SmilanskySolomon-2021}.

The case $r_\alpha\in\Q$, namely the \emph{commensurable} case, is rather different, and is the new result contained in this contribution. 
The first step is to observe that, by introducing labelled prototiles and modifying the substitution rule in a suitable way, one can construct a standard primitive substitution that ``covers'' the commensurable multiscale substitution construction, in the sense that it generates substitution tilings that are geometrically identical to the corresponding $\alpha$-Kakutani tilings described above.
 Uniform spreadness of Delone sets associated with primitive substitution tilings can be approached using Solomon's criterion \cite{Solomon-2014}. This criterion, which is a careful application of the aforementioned Laczkovich criterion, is given in terms of a comparison between the leading Perron-Frobenius eigenvalue of the associated substitution matrix, and the subsequent eigenvalues. In the special case considered here, this becomes a question of classification of certain \emph{Pisot–Vijayaraghavan (PV) polynomials} \cite{DubickasJankauskas-2014}, from which the special values that appear in \eqref{eq:special values} can be extracted.

We complete our discussion by showing that such scarcity of uniformly spread constructions does not necessarily carry over when moving away from $\alpha$-Kakutani tilings to more complicated multiscale substitution tilings. For example, when considering $(\alpha,\beta)$-Kakutani tilings, generated by a substitution rule on $I$, now split into three tiles of lengths $\alpha,\beta,1-\alpha-\beta$, we present a countable set of parameters $(\alpha,\beta)$ for which the associated Delone sets are uniformly spread.

\section{\texorpdfstring{$\alpha$}{alpha}-Kakutani tilings}

The $\alpha$-Kakutani tilings form a special family of multiscale substitution tilings of $\R$. 
The construction is inherently dynamical, being defined by the repeated application of a substitution rule. 
In this section we will introduce these tilings and state some relevant properties, for a detailed discussion of the general construction see \cite{SmilanskySolomon-2021}.

\begin{definition}\label{def:Kakutani_rule}
  Fix $\alpha\in(0,1)$ and let $I$ denote a closed interval of unit length in $\R$. The \emph{$\alpha$-Kakutani multiscale substitution rule} is a map on $I$ that returns a patch supported on $I$ that consists of two rescaled copies of $I$, a copy of $\alpha I$ to the left and a copy of $(1-\alpha)I$ to the right, see Figure \ref{Fig:Substitutionrule}. We will refer to $I$ as a \emph{prototile}. 
\begin{figure}[ht!]
	\includegraphics[scale=0.6]{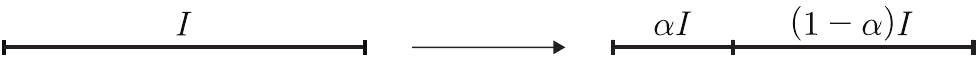}\caption{The $\alpha$-Kakutani substitution rule on $I$.}\label{Fig:Substitutionrule}
\end{figure}
We say that the $\alpha$-Kakutani substitution rule is \emph{incommensurable} if $r_\alpha\notin\Q$, and \emph{commensurable} if $r_\alpha\in\Q$, with $r_\alpha$ as defined in \eqref{eq:special values}. 
\end{definition}

We use the substitution rule to construct patches of tiles in $\R$ via the following procedure. First, position $I$ on the real line so that the origin is an interior point, and set $F_0(I)=I$ to be the patch consisting of the single tile $I$. As $t$ increases, inflate the patch by $e^t$, and substitute any tile of length greater than $1$ according to the $\alpha$-Kakutani substitution rule to define $F_t(I)$. An equivalent definition of the patch $F_t(I)$ would be to iteratively substitute $e^tI$ and all subsequent tiles of lengths greater than $1$ until they are all of at most unit length. This defines a \emph{substitution semi-flow} $F_t$. The \emph{generating family of patches} $\PP_\alpha$ is the set of all patches arising in this way, namely
$$
\PP_\alpha:=\{F_t(I)\mid t\in\R_{\geq 0}, I\text{ an interval of unit length with the origin in its interior}\}.
$$

 We define tilings of $\R$ as limits of sequences in $\PP_\alpha$ in the following way. First, identify each patch with the point set consisting of its tile boundaries, which is a closed subset of $\R$. We will use the following metric, first introduced in  \cite{Chabauty-1950}.  

\begin{definition}\label{def:CFmetric}
    Given two closed subsets $A_0,A_1\subset\R$, their \emph{Chabauty-Fell distance} $D(A_0,A_1)$ is the infimum of all $0<\varepsilon<1$ for which 
\begin{equation*}\label{eq:metric}
    A_i\cap(-1/\varepsilon,1/\varepsilon)\text{ is contained in the }\varepsilon\text{ neighbourhood of }A_{1-i}
\end{equation*}
 for both $i=0,1$, and $D(A_0,A_1)=1$ if there is no such $\varepsilon$. 
\end{definition}
Equipped with this metric, the space of closed subsets of $\R$ is compact. Simple proofs that $D$ is indeed a metric and that the space of closed subsets is compact with respect to $D$ can be found in \cite[Appendix A]{SmilanskySolomon1-2022} and \cite[\S 2]{SmilanskySolomon-2025}, respectively. The topology induced by $D$ is the Chabauty--Fell topology, also commonly referred to in the tiling literature as the local rubber topology. For a closed subset $A\subset\R$, its translation hull is
$$
    X_A=\overline{\{A+x:x\in\R\}},
$$
where the closure is taken with respect to $D$. The space $X_A$ is compact, and together with the translation action
of $\R$ it forms a topological dynamical system. When $A$ is the
boundary set of a tiling or a Delone set, we refer to $X_A$ as the corresponding space of tilings or Delone sets. See
\cite{Robinson-2004} for background on tiling dynamical systems.

\begin{definition}\label{def:Kakutani_tilings}
  An \emph{$\alpha$-Kakutani tiling} is a tiling of $\R$ whose every patch is a limit of translated sub-patches of elements in $\PP_\alpha$, where patches and tilings are identified with their set of tile-boundaries and limits are taken with respect to $D$, the Chabauty-Fell metric on closed subsets of $\R$.
\end{definition}

By construction, these tilings consist of tiles which are all intervals of lengths between $\min\{\alpha,1-\alpha\}$ and $1$. In particular, they give rise to associated Delone sets which we will denote $\Lambda_\alpha$, and are the central object of study of this paper.

\begin{figure}[ht!]
	\includegraphics[scale=0.4]{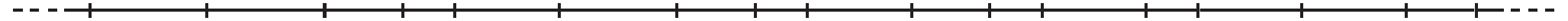}\caption{A patch of an $\alpha$-Kakutani tiling of $\R$, with $\alpha=1/3$.}\label{Fig:tiling}
\end{figure}
\begin{remark}
    The $\alpha$-Kakutani substitution rule and the procedure defining the substitution flow resemble the construction of the sequences of partitions considered by Kakutani in \cite{Kakutani-1976}, hence the name. Starting with the unit interval, in each step all sub-intervals of maximal length are split according to a fixed ratio $\alpha\in (0,1)$. A related tiling construction for the special case $\alpha=1/3$
also appears in \cite[Appendix~A.5]{FrankSadun-2014}.
\end{remark}

\begin{definition}
  Let $\alpha\in(0,1)$ and consider the $\alpha$-Kakutani substitution rule. The \emph{associated graph} $G_\alpha$ consists of a single vertex and two directed loops of lengths $\log(1/ \alpha)$ and $\log(1/(1-\alpha))$. The vertex corresponds to the tile $I$ and the edges correspond to the two rescaled copies of $I$.
  
\begin{figure}[ht!]
	\includegraphics[scale=0.6]{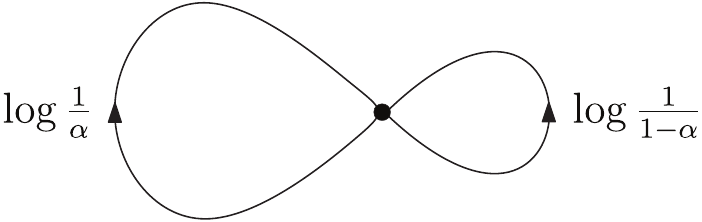}\caption{The graph associated with the $\alpha$-Kakutani substitution rule.}\label{Fig:graph}
\end{figure}
\end{definition}
Note that $r_\alpha$ is the ratio of the lengths of the two loops of
$G_\alpha$. In particular, the substitution rule is commensurable precisely when these two loop lengths are rationally related. This is a special case of the graphs associated with general multiscale substitution rules, which play a central role in the study of multiscale substitution constructions. The following result is a corollary of \cite[Proposition 2.12]{SmilanskySolomon-2021}, where the correspondence between tiles in generating patches and directed walks on the associated graph is described in detail. 
\begin{prop}\label{prop:tiles and walks}
   Let $t\in\R_{\geq 0}$. Tiles in $F_t(I)$ are in one-to-one correspondence with directed walks of length $t$ on $G_\alpha$ that originate at the single vertex. 
   If a tile $T$ in $F_t(I)$ is of length $|T|$ then the (directed) distance between the termination point of the corresponding walk on $G_\alpha$ and the vertex is $\log(1/|T|)$. In particular, the number of tiles in $F_t(I)$ is equal to the number of directed walks of length $t$ originating at the vertex of $G_\alpha$. 
\end{prop}
A directed graph is \emph{incommensurable} if it contains two closed paths of lengths $a\notin b\Q$. It is easy to see that an $\alpha$-Kakutani substitution rule is incommensurable if and only if its associated graph $G_\alpha$ is such. A multiscale substitution rule is \emph{irreducible} if $G_\alpha$ is strongly connected, which is clearly the case for $\alpha$-Kakutani substitution rules.

\section{Commensurable tilings and primitive substitutions}\label{sec: commensurable}

Commensurable multiscale substitution tilings behave rather differently than incommensurable tilings. For example, unlike incommensurable tilings, they consist of tiles of finitely many lengths. Although the translation action is minimal in both cases, commensurable tilings exhibit stronger recurrence properties such as linear repetitivity, as follows from Proposition \ref{prop:primitive cover} combined with \cite{Solomyak-1998}, compare with \cite[Theorem 1.5]{SmilanskySolomon2-2022}. We begin by showing that Delone sets associated with commensurable $\alpha$-Kakutani tilings can also be associated with primitive substitution tilings.
Such tilings, which are defined by repeated applications of a fixed scale substitution rule on an initial set of tiles, are well-studied and include such famous examples as the Penrose  and the pinwheel tilings of the plane, \cite{Penrose-1979,Radin-1994}, respectively.

\subsection{Primitive substitution tilings}
  For general definitions see \cite{BaakeGrimm-2013}. The following is tailored to our setup.
\begin{definition}
    Let $\xi>1$ and let $A=\{T_1,\ldots,T_k\}$ be a set of labelled closed intervals in $\R$, called \emph{prototiles}, each labelled by its index. A \emph{(fixed scale) substitution rule} on  $A$ with inflation constant $\xi>1$ is a map $\rho$ on $A$, so that $\rho(T_j)$ is a patch supported on $T_j$ that consists of copies of labelled rescaled prototiles in $\xi^{-1}A$, for every $T_j\in A$. 
\end{definition}
In parallel with the procedure defining $\alpha$-Kakutani tilings, given a substitution rule $\rho$ we define the \emph{generating (labelled) patches}
    $$
    \PP_\rho=\{(\xi\rho)^\ell(T)\mid \ell\in\Z_{\geq 0}, T\in A\}.
    $$
A \emph{substitution tiling} generated by $\rho$ is a tiling of $\R$ whose every patch is a translated sub-patch of an element in $\PP_\rho$, compare with Definition \ref{def:Kakutani_tilings}. By construction, all tiles in such tilings are translated copies of the prototiles in $A$. 

\begin{definition}\label{def:submatrix}
   The \emph{substitution matrix} $M_\rho=(a_{ij})$ is a $k\times k$ integer matrix with entries $a_{ij}$ given by the number of copies of $\xi^{-1}T_i$ in $\rho(T_j)$. The rule $\rho$ and the substitution tilings it generates are \emph{primitive} if $M_\rho$ is a primitive matrix, that is, if there exists $\ell\in\N$ so that $M_\rho^\ell$ is strictly positive.  
\end{definition}
The substitution matrix plays a role similar to that of the associated graph. Indeed, for every $\ell\geq 0$ and $T_j\in A$, the number of copies of $T_i$ in the generating patch $(\xi\rho)^\ell(T_j)$ is $(M_\rho^\ell)_{ij}$, compare with Proposition \ref{prop:tiles and walks}. Note that if $\rho$ is primitive, then in the definition of a substitution tiling it is enough to look only at generating patches of the form $(\xi\rho)^\ell(T_j)$, for a specific choice of prototile $T_j$.

\subsection{Commensurable tilings as primitive substitution tilings}\label{subsec:construction of rho_alpha} Recall our assumption that $\alpha=\min\{\alpha,1-\alpha\}\in(0,1/2]$, and consider a commensurable $\alpha$-Kakutani rule. Let $n\geq m \geq 1$ be coprime integers so that 
\begin{equation}\label{eq:rational ratio}
   \frac{n}{m}=r_\alpha=\frac{\log\alpha}{\log(1-\alpha)}, 
\end{equation}
or equivalently, $\alpha^m=(1-\alpha)^n$. If $n=m$ then $\alpha=1/2$, and any $\alpha$-Kakutani tiling is a translation of a lattice tiling in $\R$, with all associated Delone sets uniformly spread. Assume from now on that $n> m$, and consider the associated graph $G_\alpha$. By \eqref{eq:rational ratio} we have 
$$
g_\alpha :=\frac{1}{n}\log\frac{1}{\alpha}=\frac{1}{m}\log\frac{1}{1-\alpha}.
$$
Then by adding vertices to $G_\alpha$, the loop of length $\log(1/\alpha)$ and the loop of length $\log(1/(1-\alpha))$ can be partitioned into $n$ and $m$ edges of equal length $g_\alpha$, respectively. This defines a directed graph $G'_\alpha$ with $n+m-1$ vertices and $n+m$ edges of equal length. Label the original vertex with $1$, and assign labels $2,\ldots,n$ and $n+1,\ldots,n+m-1$ to the $(n-1)+(m-1)$ new vertices that are added along the two loops.

This adjustment of $G_\alpha$ gives rise to an adjustment of the $\alpha$-Kakutani substitution rule itself. The new vertices correspond to new labelled prototiles $\{T_2,\ldots,T_{n+m-1}\}$ and the new graph $G'_\alpha$ is associated with a new substitution rule $\rho_\alpha$, with inflation constant $\xi=e^{g_\alpha}=\alpha^{-1/n}$ on a set of prototiles $\{T_1=I,T_2,\ldots,T_{n+m-1}\}$ with
$$
|T_j|=
\begin{cases}
    |I|=1&\text{ if }\,j=1\\ \alpha^{(n+1-j)/n}&\text{ if }\,j\in\{2,\ldots,n\}\\ 
    \alpha^{(n+m-j)/n}&\text{ if }\,j\in\{n+1,\ldots,n+m-1\}.
\end{cases}   
$$
Note that $\log(1/|T_j|)$ is the length of the directed path from vertex $j$ to vertex $1$ in $G'_\alpha$. The image of $\rho_\alpha$ on these prototiles can be read from the graph $G'_\alpha$. The patch $\rho_\alpha(I)$ consists of a copy of $\alpha^{1/n}T_2$ and, if $m\geq2$, a copy of
$\alpha^{1/n}T_{n+1}$, in correspondence with the two outgoing
edges from vertex $1$. If $m=1$, the second loop has no additional vertex, and the second tile is instead a copy of
$\alpha^{1/n}I$. Similarly, for the rest of the prototiles
\begin{equation}\label{eq:primsubs}
    \rho_\alpha(T_j)=
\begin{cases}
    \alpha^{1/n}I &\text{ if }\,j\in\{n,n+m-1\}\\ \alpha^{1/n}T_{j+1} &\text{ if }\,j\in\{2,\ldots,n+m-2\}\setminus\{n\},
\end{cases} 
\end{equation}
once again in correspondence with the outgoing edges from vertex $j$. Here, if we forget the labels, geometrically speaking $\rho_\alpha$ acts on $T_1=I$ in the same way as the $\alpha$-Kakutani rule, and trivially on all other prototiles, see Example \ref{ex:r_alpha=3/2} and Figure \ref{Fig:primsubsrule} below.
Since $\gcd(n,m)=1$, the graph $G'_\alpha$ is aperiodic and therefore its adjacency matrix is primitive. The transpose of this matrix is the substitution matrix $M_{\rho_\alpha}$, which is therefore also primitive, and so $\rho_\alpha$ is a primitive substitution rule. Compare this procedure with Sadun's analysis of rational generalised pinwheel tilings \cite[\S 6]{Sadun-1998}, where a similar approach is used to study statistical properties of commensurable Kakutani substitution tilings by passing to a geometrically equivalent fixed-scale substitution rule. 

\begin{example}\label{ex:r_alpha=3/2}
    Choose $\alpha\in(0,1)$ so that $r_\alpha=3/2$, that is, $n=3,m=2$. The associated graph $G_\alpha$ can be adjusted by adding $(3-1)+(2-1)$ new vertices along the two loops to form a directed graph $G'_\alpha$ with a total of four vertices: vertex $1$ is the original vertex of $G_\alpha$, vertices $2$ and $3$ are added along the loop of length $\log(1/\alpha)$ and a vertex $4$ is added along the other loop, as illustrated in Figure \ref{Fig:commgraph}. 

    \begin{figure}[ht!]
	\includegraphics[scale=0.6]{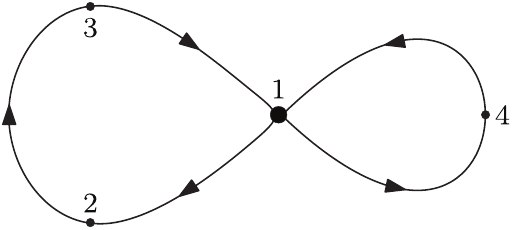}\caption{The graph $G'_\alpha$ with $r_\alpha =\frac{3}{2}$ has $4=1+(3-1)+(2-1)$ vertices and edges of equal length $g_\alpha=\frac{1}{3}\log\frac{1}{\alpha}$.}\label{Fig:commgraph}
\end{figure}
    
All edges in $G'_\alpha$ are of length $g_\alpha=\frac{1}{3}\log\frac{1}{\alpha}$. The graph gives rise to a substitution rule $\rho_\alpha$ on four prototiles $\{I,T_2,T_3,T_4\}$, given symbolically by
$$
I\mapsto T_2T_4,\qquad
T_2\mapsto T_3,\qquad
T_3\mapsto I,\qquad
T_4\mapsto I,
$$ 
whose substitution matrix is 
\[
M_{\rho_\alpha}=
\begin{pmatrix}
0&0&1&1\\
1&0&0&0\\
0&1&0&0\\
1&0&0&0
\end{pmatrix}.
\]
The prototiles are intervals of lengths $|I|=1$, $|T_2|=\alpha^{2/3}$ and $|T_3|=|T_4|=\alpha^{1/3}$, and the inflation constant is $\xi=e^{g_\alpha}=\alpha^{-1/3}$. Since $\gcd(2,3)=1$, the rule $\rho_\alpha$ is primitive. The substitution rule is illustrated in Figure \ref{Fig:primsubsrule}. 

\begin{figure}[ht!]
	\includegraphics[scale=0.6]{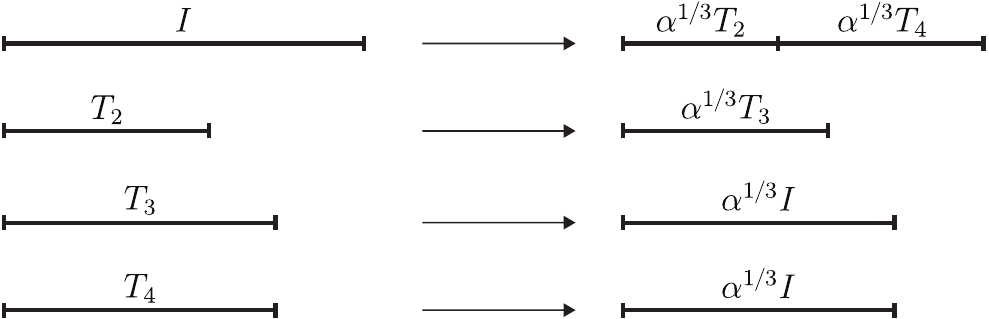}\caption{A primitive substitution rule $\rho_\alpha$ for the case $r_\alpha=\frac{3}{2}$.}\label{Fig:primsubsrule}
\end{figure}
\end{example}

\begin{prop}\label{prop:generating patches}
    Consider a commensurable $\alpha$-Kakutani substitution rule with $r_\alpha=n/m$ for $n,m\in\N$ coprime. Let $\ell\in\Z_{\geq 0}$ and put 
    \begin{equation}\label{eq:t_l}
        t_\ell:= \ell\cdot g_\alpha=\ell\cdot \frac{1}{n}\log\frac{1}{\alpha}.
    \end{equation}
    Then the (Kakutani) generating patch $F_{t_\ell}(I)$ and the (primitive substitution) generating patch $(\xi\rho_\alpha)^\ell(I)$ are geometrically the same, where $\xi=\alpha^{-1/n}$.
\end{prop}
\begin{proof}
 We prove by induction on $\ell\geq0$. For $\ell=0$, both patches contain only the single tile $I$ and so are geometrically the same. Assume the claim holds for some $\ell\geq0$. Tiles in $F_{t_\ell}(I)$ and $(\xi\rho_\alpha)^\ell(I)$ are either intervals of unit length, that is, copies of $I$, or intervals of length at most $|T_n|=|T_{n+m-1}|=\alpha^{1/n}$. By construction, the patch $F_{t_{\ell+1}}(I)$ is defined by inflating the patch $F_{t_\ell}(I)$ by $\alpha^{-1/n}$ and then substituting every tile of length strictly greater than $1$ according to the $\alpha$-Kakutani rule. The tiles that are substituted are exactly the unit intervals from $F_{t_\ell}(I)$. On the other hand, the tiles in $(\xi\rho_\alpha)^{\ell+1}(I)$ are the result of one application of $\xi\rho_\alpha$ on the tiles in $(\xi\rho_\alpha)^\ell(I)$. For copies of $I$ this amounts to an application of the $\alpha$-Kakutani rule once and inflation by $\xi=\alpha^{-1/n}$, while all other tiles in the patch are simply inflated by the same constant, and so $F_{t_{\ell+1}}(I)$ and $(\xi\rho_\alpha)^{\ell+1}(I)$ are the same.
\end{proof}

\begin{prop}\label{prop:primitive cover}
    Let $\TT$ be a commensurable $\alpha$-Kakutani tiling. Then $\TT$ is a rescaled copy of a primitive substitution tiling generated by $\rho_\alpha$, with labels forgotten. 
\end{prop}

\begin{proof}
    First, for $\ell\geq1$ and $t_\ell$ as in \eqref{eq:t_l}, if $t_{\ell-1}<t\leq t_\ell$ then $F_{t}(I)$ is a rescaled copy of $F_{t_\ell}(I)$, and so by Proposition \ref{prop:generating patches}, of $(\xi\rho_\alpha)^\ell(I)$. Indeed, this is due to the fact that all closed walks on $G_\alpha$ are of length $t_\ell$ for some $\ell$, and so by the correspondence between walks on the graph and tiles in generating patches in Proposition \ref{prop:tiles and walks}, these are the only values for which $F_t(I)$ contains a tile of unit length. 

It follows that every generating patch $F_t(I)$ is of the form $(\xi\rho_\alpha)^{\ell}(I)$, up to rescaling, consisting of intervals of finitely many lengths and with constant ratios between them. 
Since by definition $\mathcal T$ is a tiling whose every patch is a limit of translated sub-patches of generating patches, we may pass to a subsequence for which the values of $t$ modulo
$\log\xi=g_\alpha=\frac{1}{n}\log\frac{1}{\alpha}$ converge and, passing to a further subsequence if necessary, the corresponding rescaling factors converge. It follows that $\mathcal T$ is a rescaled copy of a primitive substitution tiling generated by $\rho_\alpha$, with its labels forgotten.
\end{proof}

\begin{remark}
The same graph-refinement argument applies more generally to
irreducible commensurable multiscale substitution rules in
$\R^d$. After passing to an equivalent scheme whose associated
graph has rationally related edge lengths, one may add vertices
so that all edges have a common length, obtaining a fixed-scale
substitution rule. See \cite[\S7]{Smilansky-2020} for details of
this construction.
\end{remark}

\section{Uniform spreadness of Delone sets in the real line}

Let $|U|$ denote the Lebesgue measure of a bounded measurable set $U\subset \R$. A Delone set $\Lambda\subset\R$ has \emph{asymptotic density} $d_\Lambda$ if
$$
\lim_{n\rightarrow\infty}\frac{\#(\Lambda\cap U_n)}{|U_n|}
$$
exists and is equal to $d_\Lambda$ for any \emph{van Hove sequence} $(U_n)$, that is, a sequence  of bounded measurable sets so that $|(\partial U_n)^{+\varepsilon}|/|U_n|\rightarrow0$ for all $\varepsilon>0$, where $A^{+\varepsilon}$ is the $\varepsilon$ neighbourhood of $A\subset\R$. In particular, $d_\Lambda$ is independent of the choice of van Hove
sequence. The \emph{discrepancy} of $\Lambda$ with respect to the parameter $\beta>0$ and a measurable set $U\subset\R$ is defined as 
$$
\disc{\Lambda}{\beta}{U}:=\left|\#(\Lambda \cap U)-\beta\cdot|U|\right|.
$$
A criterion established by Laczkovich \cite{Laczkovich-1992} provides a necessary and sufficient condition for a Delone set in $\R^d$ to be uniformly spread. In the case of a Delone set $\Lambda\subset \R$ with asymptotic density $d_\Lambda$, Laczkovich's criterion can be stated in the following way.
\begin{thm}\label{thm:Laczkovich}
    A Delone set $\Lambda\subset \R$ with asymptotic density $d_\Lambda$ is uniformly spread if and only if there exists a constant $C>0$ such that
    $$
    \disc{\Lambda}{d_\Lambda}{U}\leq C\cdot|(\partial U)^{+1}|
    $$
    for all bounded measurable sets $U\subset\R$.
\end{thm}

Laczkovich's criterion is a fundamental tool in the study of bounded
displacement equivalence, and was used by Solomon to establish a
criterion for uniform spreadness of primitive substitution tilings in
\cite{Solomon-2014}. The criterion is formulated in arbitrary
dimension, but we state here only the one-dimensional case that we
require. Let $\rho$ be a primitive substitution rule in $\R$ with
inflation constant $\xi$ and substitution matrix $M_\rho$. Let
$\lambda_1,\ldots,\lambda_k$ denote the eigenvalues of $M_\rho$,
ordered so that
$|\lambda_1|\geq|\lambda_2|\geq\ldots\geq|\lambda_k|$. It follows
from a simple application of the Perron-Frobenius Theorem that the
leading eigenvalue $\lambda_1$ satisfies $\lambda_1=\xi>1$, see for
example \cite[\S 5]{SmilanskySolomon-2025}.

\begin{thm}\label{thm:BD_substitution}
Let $\rho$ be a primitive substitution rule in $\R$, and let
$\Lambda$ be a Delone set associated with a tiling generated by
$\rho$. Let $\ell\geq2$ be the minimal index such that the
corresponding eigenspace satisfies
$W_{\lambda_\ell}\nsubseteq\mathbf{1}^\perp$, where
$\mathbf{1}=(1,\ldots,1)\in\R^k$. Then
\begin{enumerate}
\item
If $\absolute{\lambda_\ell}<1$ then $\Lambda$ is uniformly spread.
\item
If $\absolute{\lambda_\ell}>1$ then $\Lambda$ is not uniformly spread.
\end{enumerate}
\end{thm}

We note that the equality case $\absolute{\lambda_\ell}=1$
must be checked directly. In general, the substitution matrix alone does not determine uniform spreadness, as illustrated
by planar examples with the same matrix but different BD behaviour in \cite{FrettlohSmilanskySolomon-2021}.

\section{Proof of the main result}

Incommensurable and commensurable multiscale substitution tilings are very different in nature and their analysis requires a different set of tools. It has already been established in \cite[Theorem 8.2]{SmilanskySolomon-2021} that Delone sets associated with incommensurable multiscale substitution tilings are never uniformly spread, and so we begin with a short outline of the arguments for the special case of $\alpha$-Kakutani tilings. The analysis of the commensurable case is new and will be given in full detail.

\subsection{The incommensurable case}\label{section:incommensurable} The $\alpha$-Kakutani substitution rule is incommensurable for almost every choice of $\alpha\in(0,1)$. Although the corresponding tilings consist of tiles of infinitely many lengths, the distribution of lengths is well understood thanks to the connection with graphs as in Proposition \ref{prop:tiles and walks} and results on the distribution of walks on incommensurable graphs \cite{KiroSmilanskySmilansky-2020}. See \cite{Smilansky-2022} for explicit tile-counting formulas. These formulas, combined with the uniform frequency results in \cite[\S 9]{SmilanskySolomon-2021}, imply that the asymptotic density of any Delone set $\Lambda_\alpha$ associated with incommensurable $\alpha$-Kakutani tilings exists and is equal to
$$
d_{\Lambda_\alpha}=\frac{1}{-\alpha\log\alpha-(1-\alpha)\log(1-\alpha)}.
$$

The discrepancy can be large, however, and by \cite[\S 6]{SmilanskySolomon-2025} there exist arbitrarily large intervals $U$ for which 
\begin{equation}\label{eq:incomm_discrepancy}
    \disc{\Lambda_\alpha}{d_{\Lambda_\alpha}}{U}=\Omega\left(\frac{|U|}{\log(|U|)}\right).
\end{equation}
By Laczkovich's criterion as stated in Theorem \ref{thm:Laczkovich}, if $\Lambda_\alpha$ is uniformly spread then there exists a constant $C>0$ so that for any interval $U$ the discrepancy is bounded by $\disc{\Lambda_\alpha}{d_{\Lambda_\alpha}}{U}\leq 2C$. By \eqref{eq:incomm_discrepancy} this is clearly not the case, and so $\Lambda_\alpha$ is not uniformly spread. In fact, for any fixed incommensurable $\alpha$-Kakutani rule there are continuously many associated Delone sets that are pairwise BD non-equivalent, see \cite{SmilanskySolomon1-2022}.

\subsection{The commensurable case}
The $\alpha$-Kakutani substitution rule is commensurable for countably many choices of $\alpha\in(0,1)$. For such a choice of $\alpha$, let $\rho_\alpha$ be the substitution rule constructed in \S \ref{subsec:construction of rho_alpha}. By Proposition \ref{prop:primitive cover}, $\alpha$-Kakutani tilings can be viewed as primitive substitution tilings generated by $\rho_\alpha$. The case $r_\alpha=1$ is trivial, and we let  $n> m\geq 1$ denote coprime integers for which $r_\alpha=\frac{n}{m}$. Let $M_\alpha:=M_{\rho_\alpha}\in M_{n+m-1}(\Z)$ denote the associated substitution matrix as in Definition \ref{def:submatrix}. The matrix entries $a_{ij}$ can be read from the substitution rule $\rho_\alpha$ described in \eqref{eq:primsubs}. Indeed, we have
\begin{equation*}
\begin{cases}
a_{j+1,j}=1
&\text{if }j\in\{1,\ldots,n-1\}\sqcup
\{n+1,\ldots,n+m-2\},\\
a_{1,j}=1
&\text{if }j=n\text{ or }j=n+m-1,\\
a_{n+1,1}=1
&\text{if }m\geq2,\\
a_{1,1}=1
&\text{if }m=1,\\
a_{ij}=0
&\text{otherwise}.
\end{cases}
\end{equation*}
Thus, when $m\geq2$, the substitution matrix is of the form 
\begin{equation}\label{eq:substitution matrix}
M_{\alpha} = 
\begin{pmatrix}
0      & \cdots & 0 & 1 & 0 & \cdots & 0      & 1      \\
1      & \ddots &   &   &   &        &        & 0      \\
0      & \ddots &   &   &   &        &        & \vdots \\
       &        & 1 &   &   &        &        &        \\
1      &        &   & 0 &   &        &        &        \\
0      &        &   &   & 1 &        &        &        \\
\vdots &        &   &   &   & \ddots & \ddots &        \\
0      & \cdots &   &   &   &        & 1      &  0  
\end{pmatrix},
\end{equation}
and when $m=1$, the entry $a_{n+1,1}=1$ is replaced by
$a_{1,1}=1$.

In either case, the characteristic polynomial $p_\alpha(x)$ is easily derived by direct computation, for example, by an expansion along the first row. We deduce that
\begin{equation}\label{eq:charpoly}
p_\alpha(x)=\det(xI-M_{\alpha})=x^{n+m-1}-x^{m-1}-x^{n-1}=x^{m-1}f_\alpha(x)
\end{equation}
with $f_\alpha(x):=x^n-x^{n-m}-1$. Then the non-zero eigenvalues of the substitution matrix are precisely the roots of $f_\alpha(x)$. 
Vieta's formula implies that the product of the roots of $f_\alpha(x)$ is of modulus $1$, but since $\lambda_1=\xi>1$ is a root (it is the leading eigenvalue of $M_\alpha$, as mentioned above), there must be at least one additional non-zero eigenvalue of distinct modulus. The following results establish properties of eigenvalues of $M_\alpha$, and particularly of $\lambda_2$, which we require for our proof. 

\begin{prop}\label{prop:eigenspacesnotperp}
    Let $0\neq\lambda\in\C$ be a non-zero eigenvalue of $M_{\alpha}$. Then its eigenspace satisfies $W_\lambda\nsubseteq {\bf 1}^\perp$.
\end{prop}

\begin{proof}
    Assume otherwise, that $W_\lambda\subseteq {\bf 1}^\perp$, and let ${\bf 0}\neq v=(v_1,\ldots,v_{n+m-1})\in W_\lambda$. By assumption, $0={\bf 1}^Tv=\sum v_j$. Since the eigenspace $W_\lambda$ is $M_\alpha$-invariant, we deduce that ${\bf 1}^TM_\alpha v=0$. Note that ${\bf 1}^TM_\alpha=(2,1,\ldots,1)$, the vector of sum of column entries in $M_\alpha$. Combining the above we deduce that
    $$
    0={\bf 1}^TM_\alpha v=\sum v_j+v_1=v_1.
    $$
    But $M_\alpha v=\lambda v$, and so by observing the implied entry-wise equations we see that 
    $$
    v_1=
    \begin{cases}
        \lambda^{j-1} v_j&\text{ if }\,j\in\{2,\ldots,n\}\\ 
        \lambda^{j-n} v_j&\text{ if }\,j\in\{n+1,\ldots,n+m-1\}.
    \end{cases}
    $$
    Since $\lambda\neq0$ then $v={\bf 0}$, in contradiction with our initial assumption.
\end{proof}

\begin{prop}\label{prop:equalitycase}
Assume $|\lambda_2|=1$. Then $n=5$ and $m=1$, namely $r_\alpha=5$.
\end{prop}

Before the proof, let us recall the definition of PV numbers. 

\begin{definition}
    An algebraic integer $\lambda>1$ is a \emph{Pisot–Vijayaraghavan (PV) number} if all of its other conjugates are smaller than $1$ in modulus. An irreducible polynomial in $\Z[x]$ is a \emph{Pisot–Vijayaraghavan (PV) polynomial} if it is the minimal polynomial of a PV number. 
\end{definition}

\begin{proof}
Let $z\in\C$ be a root of $f_\alpha(x)=x^n-x^{n-m}-1$ with $|z|=1$. Solving $|z^{n-m}+1|=|z^n|=1$ we deduce that 
\begin{equation}\label{eq:z^n_z^n-m}
    z^{n-m}=e^{\pm2\pi i/3},
        \qquad
        z^n=e^{\pm\pi i/3},
\end{equation}
with the same choice of sign. Then $z^{6n}=z^{6(n-m)}=1$. Since $\gcd(n,n-m)=1$, Bezout's
identity gives $z^6=1$, and since $z^n$ is a primitive sixth
root of unity, so is $z$. Comparing exponents in \eqref{eq:z^n_z^n-m}, we see that $(n,m)\equiv(1,5)$ or $(n,m)\equiv(5,1)$ modulo $6$.

Suppose now that $|\lambda_2|=1$ and recall that $\lambda_1=\xi>1$ is also a root of $f_\alpha$. All other roots, including the other conjugates of $\xi$, have modulus at most $1$. These conjugates cannot be the primitive sixth roots found above, whose minimal polynomial is $x^2-x+1$, and are therefore strictly smaller than $1$ in modulus. Thus $\xi$ is a PV number.

Let $\psi$ denote the plastic ratio, which is the smallest PV number \cite{Siegel-1944} and the leading root of $x^3-x-1$. If $(n,m)\neq (5,1)$, then $n>5$ and so
\begin{equation*}
    \begin{split}
 f_\alpha(\psi)
 &=\psi^n-\psi^{n-m}-1\geq\psi^n-\psi^{n-1}-1\\&
  >\psi^5-\psi^4-1=(\psi^2-\psi+1)(\psi^3-\psi-1)=0.
\end{split}
\end{equation*}
Since $f_\alpha$ is strictly increasing on $[1,\infty)$, this
implies that its leading root satisfies $\xi<\psi$, a
contradiction to the minimality of $\psi$. Hence $(n,m)=(5,1)$.
\end{proof}

\begin{prop}\label{prop:r-five}
Assume $r_\alpha=5$. Then every Delone set associated with an
$\alpha$-Kakutani tiling is not uniformly spread.
\end{prop}

\begin{proof}
For $(n,m)=(5,1)$, the substitution rule is given symbolically by 
\begin{equation*}
  I\longmapsto T_2I,\qquad
  T_2\longmapsto T_3,\qquad
  T_3\longmapsto T_4,\qquad
  T_4\longmapsto T_5,\qquad
  T_5\longmapsto I.
\end{equation*}
For $k\geq0$, let
\begin{equation*}
    W_k=(\xi\rho_\alpha)^k(I)
\end{equation*}
be the patch obtained by applying the inflation--substitution map
$k$ times to $I$, and let $N_k$ be its number of tiles. We call
$W_k$ the order-$k$ supertile of type $I$, and identify it with
the word of its tile labels when writing concatenations. Since $I$
has unit length, the support of $W_k$ has length $\xi^k$. 

In our case we have
\begin{equation*}
    f_\alpha(x)=x^5-x^4-1=(x^2-x+1)(x^3-x-1).
\end{equation*}
It follows that the inflation factor $\xi>1$ satisfies $\xi^3=\xi+1$, and so for all $k\geq 0$
\begin{equation}\label{eq:supertile lengths}
    |\operatorname{supp}(W_{k+3})|
 =
 |\operatorname{supp}(W_{k+1})|
 +|\operatorname{supp}(W_{k})|,
\end{equation}
meaning that $W_{k+3}$ has the
same length as $W_{k+1}$ and $W_k$ taken together. We show that,
at every sixth level, it nevertheless contains one additional tile,
and that these excesses can be accumulated inside legal patches.

The first few supertiles are
\begin{equation*}
    \begin{aligned}
 W_0&=I,          & W_1&=T_2I,       & W_2&=T_3T_2I,\\
 W_3&=T_4T_3T_2I,       & W_4&=T_5T_4T_3T_2I,    &
 W_5&=IT_5T_4T_3T_2I=W_0W_4.
\end{aligned}
\end{equation*}
Applying the substitution $k$ additional times to the last identity
gives
\begin{equation}\label{eq:recurrence}
    W_{k+5}=W_kW_{k+4},
        \qquad
        N_{k+5}=N_k+N_{k+4}.
\end{equation}
At level $k=0$, the comparison suggested by
\eqref{eq:supertile lengths} already has an excess of one tile:
$$
N_3-N_1-N_0=4-2-1=1.
$$
The same excess reappears every six levels. Indeed, using
\eqref{eq:recurrence}, we see that
$$
\begin{aligned}
 N_{k+9}
 &=N_{k+4}+N_{k+8}
  =N_{k+4}+N_{k+3}+N_{k+7},\\
 N_{k+6}
 &=N_{k+1}+N_{k+5}
  =N_{k+1}+N_k+N_{k+4},
\end{aligned}
$$
and hence
$$
 N_{k+9}-N_{k+7}-N_{k+6}=N_{k+3}-N_{k+1}-N_k.
$$
It follows that 
\begin{equation}\label{eq:excess tiles}
    N_{6j+3}=N_{6j+1}+N_{6j}+1
\end{equation}
for all $j\geq0$, which means that at every sixth level
$W_{6j+3}$ has the same length as $W_{6j+1}$ and $W_{6j}$
taken together, but contains one additional tile.

It remains to accumulate these one-tile excesses over arbitrarily
many levels inside legal patches. For $J\geq0$ and $s=0,1,3$, set
$$
        P_{J,s}=W_{6J+s}W_{6(J-1)+s}\cdots W_{6+s}W_s.
$$
These patches are legal. Indeed, \eqref{eq:recurrence} gives
$$
 W_{k+10}=W_{k+5}W_{k+9}=W_kW_{k+4}W_{k+9},
$$
and it follows inductively that $P_{J,s}$ is a prefix of
$W_{6J+s+10}$. Summing \eqref{eq:supertile lengths} and \eqref{eq:excess tiles} for $j=0,\ldots,J$ we have
$$
\begin{aligned}
 |\operatorname{supp}(P_{J,3})|
 &=
 |\operatorname{supp}(P_{J,1})|
 +|\operatorname{supp}(P_{J,0})|,\\
 \#P_{J,3}
 &=
 \#P_{J,1}+\#P_{J,0}+J+1.               
\end{aligned}
$$
By primitivity, these legal patches occur in every generated
tiling. Let $\Lambda$ be any associated Delone set, and let
$U_{J,s}$ be the support interval of an occurrence of
$P_{J,s}$. Since $\Lambda$ contains one point corresponding
to each tile,
$$
        \#(\Lambda\cap U_{J,s})
        =\#P_{J,s}+O(1),
$$
where the error is uniformly bounded and comes only from the two
ends of $U_{J,s}$.

If $\Lambda$ were uniformly spread, with density $d$, its
discrepancy on intervals would be uniformly bounded. Subtracting
$d$ times the first identity above from the second therefore gives
$$
        J+1
        \leq
        \sum_{s\in\{0,1,3\}}
        \operatorname{disc}_{\Lambda}(d,U_{J,s})
        +O(1),
$$
which is impossible as $J\to\infty$. Hence $\Lambda$ is not
uniformly spread.
\end{proof}

\begin{prop}\label{prop:PVpolynomials}
    The second eigenvalue of $M_\alpha$ satisfies $|\lambda_2|<1$ if and only if 
    \begin{equation}\label{eq:PVpoly}
    f_\alpha (x)=x^2-x-1, x^3-x-1, x^3-x^2-1 \text{ or } x^4-x^3-1,
\end{equation}
where $f_\alpha(x)=x^n-x^{n-m}-1$ as defined in \eqref{eq:charpoly}.
\end{prop}

\begin{proof}
    Recall that $\lambda_1=\xi>1$. If $|\lambda_2|<1$ then none of the roots of $f_\alpha(x)$ have modulus $1$. By a result of Tverberg \cite[\S 5]{Tverberg-1960}, in this case $f_\alpha(x)$ is irreducible. It is therefore the minimal polynomial of $\lambda_1$ and is thus a PV polynomial. By a classification of PV polynomials due to Dubickas and Jankauskas \cite[Theorem 1.1]{DubickasJankauskas-2014}, the PV polynomials of the form $x^n-x^{n-m}-1$ are precisely those listed in \eqref{eq:PVpoly}. 
\end{proof}

\begin{proof}[Proof of Theorem \ref{thm:main}]
As already discussed above, if $\alpha=1/2$ then any $\alpha$-Kakutani tiling is a translation of a lattice tiling in $\R$, and so any associated Delone set is uniformly spread. If $r_\alpha\notin\Q$ then $\alpha$-Kakutani tilings are incommensurable, and by the discussion in \S \ref{section:incommensurable} any associated Delone set is not uniformly spread. 

Assume now that $r_\alpha=n/m$ with $n>m\geq1$ coprime integers. By Proposition \ref{prop:primitive cover}, it is enough to consider
the primitive substitution rule $\rho_\alpha$, with substitution
matrix $M_\alpha$ as in \eqref{eq:substitution matrix}. We may therefore apply Solomon's
criterion for uniform spreadness, stated above as
Theorem~\ref{thm:BD_substitution}. By
Proposition \ref{prop:eigenspacesnotperp} and the preceding discussion, the eigenvalue $\lambda_2$ of $M_\alpha$ is non-zero and satisfies $W_{\lambda_2}\nsubseteq{\bf 1}^\perp$, and so $\ell=2$. If $|\lambda_2|>1$, Theorem~\ref{thm:BD_substitution} shows that $\Lambda_\alpha$ is
not uniformly spread. If $|\lambda_2|=1$, Proposition~\ref{prop:equalitycase} gives $n=5$ and $m=1$, and Proposition~\ref{prop:r-five} shows
that $\Lambda_\alpha$ is again not uniformly spread. It follows
that $\Lambda_\alpha$ is uniformly spread if and only if
$|\lambda_2|<1$. By Proposition \ref{prop:PVpolynomials}, for this to hold $f_\alpha(x)$ must be one of the four polynomials listed in \eqref{eq:PVpoly}. We thus conclude that the values of $\alpha\in(0,1/2)$ for which  $\Lambda_\alpha$ is uniformly spread are given by
\begin{equation}\label{eq:PVvalues}
   r_\alpha=\frac{n}{m}\in\left\{\frac32,\frac21,\frac31,\frac41\right\}, 
\end{equation}
and the proof is complete. 
\end{proof}

\section{Concluding remarks}\label{sec:remarks}
We now take a closer look at the four special values of $\alpha$ listed in \eqref{eq:PVvalues}. Recall that
$$
r_\alpha=\frac{n}{m}=\frac{\log(\alpha)}{\log(1-\alpha)},\quad f_\alpha(x)=x^n-x^{n-m}-1.
$$

\begin{enumerate}
    \item $n=3,m=2$. Then $f_\alpha(x)=x^3-x-1$ is the minimal polynomial of the \emph{plastic ratio}, the smallest PV number \cite{Siegel-1944}. The value of $\alpha$ satisfies $$\alpha^2=(1-\alpha)^3,\quad\alpha\approx0.43016.$$  
     
    \item $n=2,m=1$. Then $f_\alpha(x)=x^2-x-1$ is the minimal polynomial of the \emph{golden ratio} $\phi$. The value of $\alpha$ satisfies $$\alpha=(1-\alpha)^2,\quad\alpha=\frac{1}{\phi^2}\approx0.38196.$$
    
    \item $n=3,m=1$. Then $f_\alpha(x)=x^3-x^2-1$ is the minimal polynomial of the \emph{supergolden ratio}. The value of $\alpha$ satisfies $$\alpha=(1-\alpha)^3,\quad\alpha\approx0.31767.$$
    
    \item $n=4,m=1$. Then $f_\alpha(x)=x^4-x^3-1$ is the minimal polynomial of the second smallest PV number \cite{Siegel-1944}. The value of $\alpha$ satisfies $$\alpha=(1-\alpha)^4,\quad\alpha\approx0.27551.$$
\end{enumerate}

These four PV numbers are well-studied and are of interest in number theory, combinatorics and geometry. We refer to the individual Wikipedia entries of the plastic, golden and supergolden ratios for a friendly introduction and illuminating illustrations.

 \subsection{Beyond \texorpdfstring{$\alpha$}{alpha}-Kakutani tilings} Let us now consider other families of multiscale substitution tilings parametrised by continuous parameters. One such family of tilings of $\R$ is defined similarly to $\alpha$-Kakutani tilings, but with respect to a substitution rule that splits the unit interval into three intervals of lengths $\alpha,\beta,1-\alpha-\beta$.  An analysis similar to the $\alpha$-Kakutani case leads us to consider the countably many commensurable substitution rules with pairs $\alpha,\beta$ for which there exist integers $n\geq m\geq k\geq 1$ with $\gcd(n,m,k)=1$ satisfying
$$
    \alpha^m=\beta^n,\qquad
    \alpha^k=(1-\alpha-\beta)^n.
$$
We find a primitive substitution rule that generates such tilings, with non-zero eigenvalues of the substitution matrix given as the roots of the polynomial 
$$
f_{\alpha,\beta}(x)=x^n-x^{n-m}-x^{n-k}-1.
$$
However, unlike the $\alpha$-Kakutani case, the classification in \cite{DubickasJankauskas-2014} provides infinitely many such PV polynomials: $x^5-x^4-x^2-1$ and the infinite families given by  
$$
x^d-2x^{d-1}-1 \text{ for }d\geq1 \text{ and } x^d-x^{d-1}-x^{d-2}-1 \text{ for odd }d\geq3,
$$
defining infinitely many $\alpha,\beta$ for which all associated Delone sets are uniformly spread.

\begin{figure}[ht!]
	\includegraphics[scale=0.6]{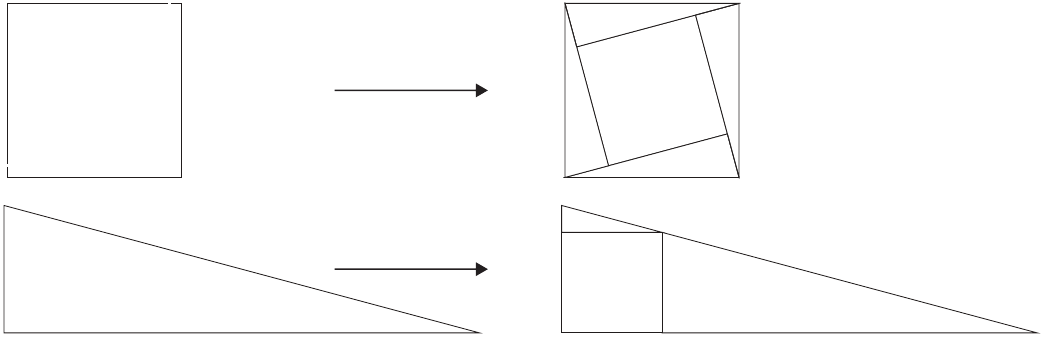}\caption{A multiscale substitution rule on a square and a right  triangle. The acute angle may be changed continuously in the range $(0,\pi/4)$.}\label{Fig:rose}
\end{figure}

We conclude our discussion by mentioning families of multiscale substitution tilings of higher dimension that continuously depend on a parameter. Examples include the aforementioned generalised pinwheel tilings \cite{Sadun-1998}, or those generated by the substitution rule illustrated in Figure \ref{Fig:rose}. 

As in the one-dimensional case, the only possible candidates for
uniformly spread Delone sets are the countably many angles that
produce commensurable tilings, and appropriate fixed-scale
substitutions can be found. However, since $d=2$, Solomon's
criterion no longer translates to a classification of PV
polynomials, and a different approach is required.

\subsection*{Acknowledgements}
I am happy to thank Alon Nishry, Yaar Solomon and Barak Weiss for helpful discussions, and the anonymous referees for their valuable comments. 
\printbibliography

%\bibliographystyle{alpha}
%\bibliography{bib_BD_Kakutani}

\end{document}